\documentclass[12pt]{article}%
\usepackage{graphicx}
\usepackage[intlimits]{amsmath}
\usepackage{latexsym}
\usepackage{amsfonts}
\usepackage{amssymb}%
\makeatletter
\newfont{\footsc}{cmcsc10 at 8truept}
\newfont{\footbf}{cmbx10 at 8truept}
\newfont{\footrm}{cmr10 at 10truept}
\newtheorem{theorem}{Theorem}

\newtheorem{lemma}[theorem]{Lemma}

\newtheorem{proposition}[theorem]{Proposition}

\newenvironment{proof}[1][Proof]{\noindent{\textbf {#1}  }}  {\hfill$\Box$\bigskip}

\begin{document}

\title{The sum of the squares of degrees: an overdue assignement }
\author{Vladimir Nikiforov\\{\small Department of Mathematical Sciences, University of Memphis, Memphis,
Tennessee, 38152}}
\maketitle

\begin{abstract}
Let $f\left(  n,m\right)  $ be the maximum of the sum of the squares of
degrees of a graph with $n$ vertices and $m$ edges. Summarizing earlier
research, we present a concise, asymptotically sharp upper bound on $f\left(
n,m\right)  $, better than the bound of de Caen for almost all $n$ and $m.$

\textbf{Keywords:}\textit{ squares of degrees, de Caen's bound, sharp bound}

\end{abstract}

\section{Introduction}

Our notation is standard (e.g., see \cite{Bol98}). Specifically, in this note,
$n$ and $m$ denote the number of vertices and edges of a graph $G$.

Few problems in combinatorics have got so many independent solutions as the
problem of finding
\[
f\left(  n,m\right)  =\max\left\{  \sum_{u\in V\left(  G\right)  }d^{2}\left(
u\right)  :v\left(  G\right)  =n,\text{ }e\left(  G\right)  =m\right\}  .
\]

The first contribution is due to B. Schwarz \cite{Sch64} who studied how to
shuffle the entries of a square nonnegative matrix $A$ in order to maximize
the sum of the entries of $A^{2}$. Later M. Katz \cite{Kat71} almost
completely solved the same problem for square $\left(  0,1\right)  $-matrices,
obtaining, in particular, an asymptotic value of $f\left(  n,m\right)  .$ The
first exact result for $f\left(  n,m\right)  ,$ found in 1978 by Ahlswede and
Katona \cite{AlKa78}, reads as: suppose $r,q,s,t$ are integers defined
uniquely by
\begin{equation}
m=\binom{r}{2}+q=\binom{n}{2}-\binom{s}{2}-t,\text{ \ \ }0\leq q<r,\text{
}0\leq t<s, \label{eq1}%
\end{equation}
and set
\begin{align}
C\left(  n,m\right)   &  =2m\left(  r-1\right)  +q\left(  q+1\right)
,\label{defc}\\
S\left(  n,m\right)   &  =\left(  n\left(  n-1\right)  -2m\right)  \left(
s-1\right)  +t\left(  t+1\right)  +4m\left(  n-1\right)  -\left(  n-1\right)
^{2}n. \label{defs}%
\end{align}
Then
\begin{equation}
f\left(  n,m\right)  =\max\left\{  C\left(  n,m\right)  ,S\left(  n,m\right)
\right\}  . \label{form}%
\end{equation}

Moreover, Ahlswede and Katona demonstrated that, if $\left\vert m-n\left(
n-1\right)  /4\right\vert <n/2,$ finding $\max\left\{  C\left(  n,m\right)
,S\left(  n,m\right)  \right\}  $ is a subtle and difficult problem; hence,
there is little hope for a simple exact expression for $f\left(  n,m\right)
.$

Almost at the same time Aharoni \cite{Aha79} completed the work of Katz for
square $\left(  0,1\right)  $-matrices. In 1987 Brualdi and Solheid
\cite{BrSo87}, adapting Aharoni's method to graphs, rediscovered (\ref{form})
and in 1996 Olpp \cite{Olp96}, apparently unaware of these achievements,
meticulously deduced (\ref{form}) from scratch.

Despite this impressive work, none of these authors came up with a concise,
albeit approximate upper bound on $f\left(  n,m\right)  .$ In contrast, de
Caen \cite{Cae98} proved that
\begin{equation}
f\left(  n,m\right)  \leq m\left(  \frac{2m}{n-1}+n-2\right)  . \label{Cae}%
\end{equation}

Denote the right-hand side of (\ref{Cae}) by $D\left(  n,m\right)  $ and note
that, for almost all $n$ and $m,$ it is considerably greater than $f\left(
n,m\right)  $ - in fact, for $m$ around $n^{2}/4$ and $n$ sufficiently large,
$D\left(  n,m\right)  >1.06f\left(  n,m\right)  .$ De Caen was aware that
$D\left(  n,m\right)  $ matches $f\left(  n,m\right)  $ poorly, but he
considered that it has \textquotedblleft... an appealingly simple
form.\textquotedblright\ He was right - his result motivated further research,
e.g., see \cite{Bey03}, \cite{Cio06}, and \cite{Das04}. Sadly enough, neither
de Caen, nor his successors refer to the work done before Olpp.

In summary: the result (\ref{form}) is exact but complicated, while de Caen's
result (\ref{Cae}) is simple but inexact.

The aim of this note is to find a concise, asymptotically sharp upper bound on
$f\left(  n,m\right)  $, better than de Caen's bound for almost all $n$ and
$m.$

We begin with the following \textquotedblleft half\textquotedblright\ result.

\begin{theorem}
\label{th1}If $m\geq n\left(  n-1\right)  /4,$ then%
\begin{equation}
m\sqrt{8m+1}-3m\leq f\left(  n,m\right)  \leq m\sqrt{8m+1}-m. \label{bo1}%
\end{equation}
Moreover, for $m<\left(  n-1\right)  \left(  n-2\right)  /2,$
\begin{equation}
m\sqrt{8m+1}-m<D\left(  n,m\right)  . \label{bo2}%
\end{equation}

\end{theorem}

This theorem is almost as good as one can get, but it holds only for half of
the range of $m$. Since
\[
f\left(  n,\frac{n\left(  n-1\right)  }{2}-m\right)  =f\left(  n,m\right)
+4\left(  n-1\right)  m-n\left(  n-1\right)  ^{2},
\]
one can produce a bound when $m<n\left(  n-1\right)  /4$ as well. We state
below a simplified complete version.

\begin{theorem}
\label{th2}Let
\[
F\left(  n,m\right)  =\left\{
\begin{array}
[c]{ll}%
\left(  2m\right)  ^{3/2}, & if\text{ }m\geq n^{2}/4\\
\left(  n^{2}-2m\right)  ^{3/2}+4mn-n^{3}, & if\text{ }m<n^{2}/4.
\end{array}
\right.
\]
Then, for all $n$ and $m,$
\begin{equation}
F\left(  n,m\right)  -4m\leq f\left(  n,m\right)  \leq F\left(  n,m\right)  .
\label{bo3}%
\end{equation}
Moreover, if $n^{3/2}<m<\binom{n}{2}-n^{3/2},$ then
\begin{equation}
F\left(  n,m\right)  <D\left(  n,m\right)  . \label{bo4}%
\end{equation}

\end{theorem}

\section{Proofs}

To begin with, note that (\ref{defc}) and (\ref{defs}) imply that
\begin{equation}
S\left(  n,m\right)  =C\left(  n,\frac{n\left(  n-1\right)  }{2}-m\right)
+4m\left(  n-1\right)  -n\left(  n-1\right)  ^{2}. \label{sc}%
\end{equation}

We need some preliminary results.

\begin{proposition}
\label{pro2}For all $n$ and $m>0,$
\begin{equation}
\left(  2m\right)  ^{3/2}-3m<m\sqrt{8m+1}-3m\leq C\left(  n,m\right)  .
\label{p1}%
\end{equation}

\end{proposition}

\begin{proof}
Let $m=\binom{r}{2}+q,$ $0\leq q<r.$ From
\[
\left(  8m\right)  ^{1/2}<\sqrt{8m+1}=\sqrt{4r\left(  r-1\right)  +8q+1}<2r+1
\]
and (\ref{defc}) we deduce that
\[
C\left(  n,m\right)  =2m\left(  r-1\right)  +q\left(  q+1\right)
\geq2m\left(  r+\frac{1}{2}\right)  -3m\geq m\sqrt{8m+1}-3m,
\]
proving (\ref{p1}) and the proposition.
\end{proof}

\begin{proposition}
\label{pr0} For every $r\geq3$
\[
\sqrt{\left(  2r-1\right)  ^{2}+8\left(  r-1\right)  }>\frac{2r^{2}+5r-2}%
{r+2}.
\]

\end{proposition}

\begin{proof}
Since%
\[
\sqrt{\left(  2r-1\right)  ^{2}+8\left(  r-1\right)  }=\sqrt{\left(
2r+1\right)  ^{2}-8},
\]
the desired inequality follows from
\begin{align*}
\left(  2r+1\right)  ^{2}-8  &  \geq\left(  2r+1\right)  ^{2}-8\frac
{2r^{2}+5r}{\left(  r+2\right)  ^{2}}\geq\left(  2r+1\right)  ^{2}%
-8\frac{\left(  2r+1\right)  \left(  r+2\right)  -2}{\left(  r+2\right)  ^{2}%
}\\
&  =\left(  2r+1\right)  ^{2}-\frac{8\left(  2r+1\right)  }{r+2}+\frac
{16}{\left(  r+2\right)  ^{2}}=\left(  \frac{2r^{2}+5r-2}{r+2}\right)  ^{2},
\end{align*}
completing the proof.
\end{proof}

\begin{lemma}
\label{pro1}For all $n$ and $m,$%
\[
C\left(  n,m\right)  \leq m\sqrt{8m+1}-m.
\]

\end{lemma}

\begin{proof}
Let $m=\binom{r}{2}+q,$ $0\leq q<r.$ In view of (\ref{eq1}) and (\ref{defc}%
),\ the required inequality is equivalent to
\[
2r\left(  r-1\right)  ^{2}+4rq+2q\left(  q-1\right)  \leq\left(  r\left(
r-1\right)  +2q\right)  \sqrt{\left(  2r-1\right)  ^{2}+8q}-r\left(
r-1\right)  -2q,
\]
and so, to%
\begin{equation}
\left(  2r-1\right)  r\left(  r-1\right)  \leq\left(  r\left(  r-1\right)
+2q\right)  \sqrt{\left(  2r-1\right)  ^{2}+8q}-4rq-2q^{2}. \label{in1}%
\end{equation}

It is immediate to check that (\ref{in1}) holds if $r=1;$ thus we shall assume
that $r\geq2$. If $q=r-1,$ then Proposition \ref{pr0} implies (\ref{in1}) by%
\begin{align*}
&  \left(  r\left(  r-1\right)  +2\left(  r-1\right)  \right)  \sqrt{\left(
2r-1\right)  ^{2}+8\left(  r-1\right)  }-4r\left(  r-1\right)  -2\left(
r-1\right)  ^{2}\\
&  =\left(  r-1\right)  \left(  \left(  r+2\right)  \sqrt{\left(  2r-1\right)
^{2}+8\left(  r-1\right)  }-6r+2\right) \\
&  >\left(  r-1\right)  \left(  2r^{2}+5r-2-6r+2\right)  =\left(  r-1\right)
r\left(  2r-1\right)  .
\end{align*}

Assume now $r\geq2,$ and $0\leq q\leq r-2.$ Then Bernoulli's inequality
implies that%
\begin{align*}
\left(  \left(  2r-1\right)  ^{2}+8q\right)  ^{3/2}  &  \geq\left(
2r-1\right)  ^{3}\left(  1+\frac{12q}{\left(  2r-1\right)  ^{2}}\right)
=\left(  2r-1\right)  ^{3}+12q\left(  2r-1\right)  ,\\
\left(  \left(  2r-1\right)  ^{2}+8q\right)  ^{1/2}  &  \leq\left(
2r-1\right)  \left(  1+\frac{4q}{\left(  2r-1\right)  ^{2}}\right)  =\left(
2r-1\right)  +\frac{4q}{\left(  2r-1\right)  },
\end{align*}
and so,
\begin{align*}
&  \left(  r\left(  r-1\right)  +2q\right)  \sqrt{\left(  2r-1\right)
^{2}+8q}-4rq-2q^{2}\\
&  =\frac{1}{4}\left(  \left(  2r-1\right)  ^{2}+8q\right)  ^{3/2}-\frac{1}%
{4}\left(  \left(  2r-1\right)  ^{2}+8q\right)  ^{1/2}-4rq-2q^{2}\\
&  >\frac{\left(  2r-1\right)  ^{3}+12q\left(  2r-1\right)  }{4}-\frac{\left(
2r-1\right)  }{4}-\frac{q}{\left(  2r-1\right)  }-4rq-2q^{2}\\
&  =\left(  2r-1\right)  r\left(  r-1\right)  +q\left(  2r-3-2q-\frac
{1}{\left(  2r-1\right)  }\right) \\
&  \geq\left(  2r-1\right)  r\left(  r-1\right)  +q\left(  2r-3-2\left(
r-2\right)  -\frac{1}{\left(  2r-1\right)  }\right)  \geq\left(  2r-1\right)
r\left(  r-1\right)  .
\end{align*}

This completes the proof of (\ref{in1}) and of Lemma \ref{pro1}.
\end{proof}

\begin{proof}
[\textbf{Proof of Theorem \ref{th1}}]The first inequality in (\ref{bo1})
follows from $C\left(  n,m\right)  \leq f\left(  n,m\right)  $ and Proposition
\ref{pro2}. To prove the second inequality in (\ref{bo1}), set first
\[
A\left(  n,m\right)  =\left(  \frac{n\left(  n-1\right)  }{2}-m\right)
\sqrt{\left(  2n-1\right)  ^{2}-8m}-\frac{n\left(  n-1\right)  }%
{2}+m+4m\left(  n-1\right)  -n\left(  n-1\right)  ^{2}%
\]
and observe that (\ref{sc}) and Lemma \ref{pro1} imply that, for all $n$ and
$m,$%
\begin{equation}
S\left(  n,m\right)  \leq A\left(  n,m\right)  . \label{i00}%
\end{equation}
We shall prove that, if $m\geq n\left(  n-1\right)  /4,$ then\
\begin{equation}
A\left(  n,m\right)  \leq m\sqrt{8m+1}. \label{i0}%
\end{equation}
Setting $x=\frac{n\left(  n-1\right)  }{2}-m,$ this is equivalent to: if
$x\leq n\left(  n-1\right)  /4,$ then
\begin{align}
&  x\sqrt{8x+1}-x-4x\left(  n-1\right)  +n\left(  n-1\right)  ^{2}\nonumber\\
&  \leq\left(  \frac{n\left(  n-1\right)  }{2}-x\right)  \sqrt{8\left(
\frac{n\left(  n-1\right)  }{2}-x\right)  +1}-\frac{n\left(  n-1\right)  }%
{2}+x. \label{i1}%
\end{align}
Setting $g\left(  x\right)  =x\sqrt{8x+1}-\left(  2n-1\right)  x,$ (\ref{i1})
is equivalent to: if $0\leq x\leq n\left(  n-1\right)  /4,$ then
\[
g\left(  x\right)  \leq g\left(  \frac{n\left(  n-1\right)  }{2}-x\right)
\]
Since,
\[
g^{\prime}\left(  x\right)  =\sqrt{8x+1}+4x\left(  8x+1\right)  ^{-1/2}%
-\left(  2n-1\right)  \geq4x\left(  8x+1\right)  ^{-1/2}>0,
\]
$g\left(  x\right)  $ increases with $x,$ and $g\left(  \frac{n\left(
n-1\right)  }{2}-x\right)  $ decreases with $x.$ Hence,
\[
g\left(  x\right)  \leq g\left(  n\left(  n-1\right)  /4\right)  \leq g\left(
\frac{n\left(  n-1\right)  }{2}-x\right)  ,
\]
proving (\ref{i1}) and (\ref{i0}). Finally, if $m\geq n\left(  n-1\right)
/4,$ then Lemma \ref{pro1}, \ref{i00}, and (\ref{i0}) imply that
\[
\max\left\{  C\left(  n,m\right)  ,S\left(  n,m\right)  \right\}  \leq
\max\left\{  m\sqrt{8m+1}-m,A\left(  n,m\right)  \right\}  =m\sqrt{8m+1}-m.
\]
This, in view of (\ref{form}), completes the proof of the second inequality in
(\ref{bo1}).

\textbf{Proof of (\ref{bo2})}

To prove (\ref{bo2}), assume that $m\sqrt{8m+1}-m\geq D\left(  n,m\right)  .$
Then
\[
\frac{2m}{n-1}+n-1\leq\sqrt{8m+1}%
\]
and so,%
\[
4m^{2}-4m\left(  n-1\right)  ^{2}+n\left(  n-1\right)  ^{2}\left(  n-2\right)
\leq0,
\]
implying that
\[
\frac{2m}{n-1}\geq n-2,
\]
a contradiction with the assumption about $m$. This completes the proof of
Theorem \ref{th1}.
\end{proof}

To simplify the proof of Theorem \ref{th2}, we need the following lemma.

\begin{lemma}
\label{pro3}For $m\leq n^{2}/4,$
\begin{equation}
S\left(  n,m\right)  \leq\left(  n^{2}-2m\right)  ^{3/2}+4mn-n^{3}.
\label{in5}%
\end{equation}

\end{lemma}

\begin{proof}
Let $\binom{n}{2}-m=\binom{s}{2}+t.$ Lemma \ref{pro1} implies that
\begin{align*}
C\left(  n,\binom{n}{2}-m\right)   &  =2\left(  \binom{n}{2}-m\right)  \left(
s-1\right)  +t\left(  t+1\right) \\
&  \leq\left(  \binom{n}{2}-m\right)  \sqrt{\left(  2n-1\right)  ^{2}%
-8m}-\binom{n}{2}+m.
\end{align*}
Hence, in view of (\ref{sc}), inequality (\ref{in5}) follows from
\begin{align*}
&  \left(  \binom{n}{2}-m\right)  \sqrt{\left(  2n-1\right)  ^{2}-8m}%
-\binom{n}{2}+m+4m\left(  n-1\right)  -\left(  n-1\right)  ^{2}n\\
&  \leq\left(  n^{2}-2m\right)  ^{3/2}+4mn-n^{3},
\end{align*}
in turn, equivalent to%
\begin{equation}
2\left(  n^{2}-2m\right)  ^{3/2}-\left(  n\left(  n-1\right)  -2m\right)
\sqrt{\left(  2n-1\right)  ^{2}-8m}+6m-3n^{2}+3n\geq2n. \label{in2}%
\end{equation}
Thus, our goal is the proof of (\ref{in2}). Note the for $n\leq3$ $,$
inequality (\ref{in2}) holds for every $m$, so we shall assume that $n\geq4.$
Let%
\[
g\left(  x\right)  =2\left(  x+n\right)  ^{3/2}-x\left(  4x+1\right)
^{1/2}-3x
\]
and observe that (\ref{in2}) is equivalent to $g\left(  n\left(  n-1\right)
-2m\right)  \geq2n.$ We first prove that $g\left(  x\right)  \ $is decreasing
for $n\left(  n-1\right)  -n^{2}/2\leq x\leq n\left(  n-1\right)  $. Indeed,
\begin{align*}
g^{\prime}\left(  x\right)   &  =3\left(  x+n\right)  ^{1/2}-\left(
4x+1\right)  ^{1/2}-2x\left(  4x+1\right)  ^{-1/2}-3\\
&  =3\left(  x+n\right)  ^{1/2}-\left(  4x+1\right)  ^{1/2}-2x\left(
4x+1\right)  ^{-1/2}-3\\
&  \leq3x^{1/2}\left(  1+\frac{n}{2x}\right)  -\frac{6x+1}{\sqrt{4x+1}}%
-3\leq3x^{1/2}\left(  1+\frac{n}{2x}\right)  -\frac{6x+1}{2x^{1/2}\left(
1+1/8x\right)  }-3\\
&  =3x^{1/2}+\frac{3n}{2x^{1/2}}-\frac{24x+4}{8x+1}x^{1/2}-3<3x^{1/2}%
+\frac{3n}{2x^{1/2}}-3x^{1/2}-3\\
&  =3\frac{n}{2x^{1/2}}-3=\frac{3}{x^{1/2}}\left(  \frac{n}{2}-\left(
\frac{n^{2}}{2}-n\right)  ^{1/2}\right)  <\frac{3}{x^{1/2}}\left(  \frac{n}%
{2}-\frac{n}{\sqrt{2}}\left(  1-\frac{1}{n}\right)  \right)  <0.
\end{align*}
Therefore,
\[
g\left(  n\left(  n-1\right)  -2m\right)  \geq g\left(  n\left(  n-1\right)
\right)  =2n^{3}-n\left(  n-1\right)  \left(  2n-1\right)  -3n\left(
n-1\right)  =2n,
\]
proving (\ref{in2}) and Lemma \ref{pro3}.
\end{proof}

\begin{proof}
[\textbf{Proof of Theorem \ref{th2}}]Our first goal is to prove the second
inequality in (\ref{bo3}). Note that the function $g\left(  x\right)
=x^{3/2}-x$ is increasing for $1/2\leq x\leq1$. Indeed, $g^{\prime}\left(
x\right)  =\frac{3}{2}x^{1/2}-1>\frac{3}{2\sqrt{2}}-1>0.$ Hence, $g\left(
1-x\right)  $ is decreasing for $1/2\leq x\leq1.$ Hence, if $1/2\leq x\leq1,$
then
\[
g\left(  x\right)  \geq g\left(  1/2\right)  \geq g\left(  1-x\right)  ;
\]
likewise, if $0\leq x\leq1/2,$ then
\[
g\left(  1-x\right)  \geq g\left(  1/2\right)  \geq g\left(  x\right)  .
\]
Therefore, setting $x=2m/n^{2},$ we see that, if $n^{2}/4\leq m\leq n\left(
n-1\right)  ,$ then%
\[
\left(  2m\right)  ^{3/2}\geq\left(  n^{2}-2m\right)  ^{3/2}+4mn-n^{3}%
\]
and, if $0\leq m\leq n^{2}/4,$ then%
\[
\left(  2m\right)  ^{3/2}\leq\left(  n^{2}-2m\right)  ^{3/2}+4mn-n^{3}.
\]
In other words,%
\[
F\left(  n,m\right)  =\max\left\{  \left(  2m\right)  ^{3/2},\left(
n^{2}-2m\right)  ^{3/2}+4mn-n^{3}\right\}  .
\]

Lemma \ref{pro1} implies that, for all $n$ and $m,$
\[
C\left(  n,m\right)  \leq m\sqrt{8m+1}-m\leq\left(  2m\right)  ^{3/2};
\]
Lemma \ref{pro3} implies that, for $m\leq n^{2}/4$,
\[
S\left(  n,m\right)  \leq\left(  n^{2}-2m\right)  ^{3/2}+4mn-n^{3},
\]
and so, in view of (\ref{form}), the second inequality in (\ref{bo3}) is proved.

\textbf{Proof of the first inequality in (\ref{bo3}) }

To prove the first inequality in (\ref{bo3}), assume first that $m<n^{2}/4;$
we shall prove that%
\[
\left(  n^{2}-2m\right)  ^{3/2}-n^{3}+4mn-4m\leq S\left(  n,m\right)  .
\]
Letting $\binom{n}{2}-m=\binom{s}{2}+t,$ in view of (\ref{defs}), this is
equivalent to%
\begin{equation}
\left(  n^{2}-2m\right)  ^{3/2}\leq\left(  n\left(  n-1\right)  -2m\right)
\left(  s-1\right)  +t\left(  t+1\right)  +2n^{2}-n, \label{in3}%
\end{equation}
Thus, our goal is to prove (\ref{in3}).

Bernoulli's inequality implies that%
\begin{align*}
\left(  n\left(  n-1\right)  -2m\right)  ^{3/2}  &  =\left(  n^{2}-2m\right)
^{3/2}\left(  1-\frac{n}{n^{2}-2m}\right)  ^{3/2}\\
&  \geq\left(  n^{2}-2m\right)  ^{3/2}\left(  1-\frac{3n}{2\left(
n^{2}-2m\right)  }\right)  =\left(  n^{2}-2m\right)  ^{3/2}-\frac{3}%
{2}n\left(  n^{2}-2m\right)  ^{1/2},
\end{align*}
and so,%
\begin{equation}
\left(  n^{2}-2m\right)  ^{3/2}\leq\left(  n\left(  n-1\right)  -2m\right)
^{3/2}+\frac{3}{2}n\sqrt{n^{2}-2m}\leq\left(  n\left(  n-1\right)  -2m\right)
^{3/2}+\frac{3\sqrt{2}}{4}n^{2}. \label{in4}%
\end{equation}
On the other hand, from
\[
n\left(  n-1\right)  -2m=s\left(  s-1\right)  +2t<s\left(  s+1\right)
\]
we see that $\sqrt{n\left(  n-1\right)  -2m}<s+1/2.$ Hence, in view of
(\ref{in4}), we have
\begin{align*}
\left(  n^{2}-2m\right)  ^{3/2}  &  \leq\left(  n\left(  n-1\right)
-2m\right)  \left(  s-1\right)  +\frac{3}{2}\left(  n\left(  n-1\right)
-2m\right)  +\frac{3\sqrt{2}}{4}n^{2}\\
&  \leq\left(  n\left(  n-1\right)  -2m\right)  \left(  s-1\right)  +\frac
{3}{2}n\left(  n-1\right)  -\frac{3n^{2}}{4}+\frac{3\sqrt{2}}{4}n^{2}\\
&  <\left(  n\left(  n-1\right)  -2m\right)  \left(  s-1\right)  +2n^{2}-n,
\end{align*}
completing the proof of (\ref{in3}). Since, by Proposition \ref{pro2}, we
have
\[
\left(  2m\right)  ^{3/2}-3m\leq C\left(  n,m\right)  ,
\]
it follows that
\[
F\left(  n,m\right)  -4m\leq\left\{
\begin{array}
[c]{ll}%
C\left(  n,m\right)  , & if\text{ }m\geq n^{2}/4\\
S\left(  n,m\right)  , & if\text{ }m<n^{2}/4.
\end{array}
\right.
\]
implying the first inequality in (\ref{bo3}).

\textbf{Proof of (\ref{bo4})}

To prove (\ref{bo4}), suppose first that $n^{2}/4\leq m<\binom{n}{2}-\left(
n-1\right)  ^{3/2};$ then we have to prove that
\begin{equation}
\left(  2m\right)  ^{3/2}<m\left(  \frac{2m}{n-1}+n-2\right)  \label{in6}%
\end{equation}
Assuming that (\ref{in6}) fails, we see that
\[
2\sqrt{2m}\geq\frac{2m}{n-1}+n-2,
\]
and so,
\[
\left(  \sqrt{\frac{2m}{n-1}}-\sqrt{n-1}\right)  ^{2}\leq1.
\]
After some algebra we obtain%
\[
2m\geq n\left(  n-1\right)  -2\left(  n-1\right)  \sqrt{n-1},
\]
a contradiction with the range of $m.$

Suppose now that $n^{3/2}<m\leq n^{2}/4.$ This implies
\[
n^{2}-2\left(  n-1\right)  ^{3/2}>n^{2}-2m>n^{2}/2,
\]
and thus, by (\ref{in6}),
\[
\left(  n^{2}-2m\right)  ^{3/2}<\left(  n^{2}-2m\right)  \left(
\frac{2\left(  n^{2}-2m\right)  }{n-1}+n-2\right)  .
\]
Hence,%
\begin{align*}
&  \left(  n^{2}-2m\right)  ^{3/2}+4mn-n^{3}\\
&  <\frac{\left(  n^{2}-2m\right)  }{2}\left(  \frac{\left(  n^{2}-2m\right)
}{n-1}+n-2\right)  +4mn-n^{3}\\
&  =\frac{n^{4}-4mn^{2}+4m^{2}}{2\left(  n-1\right)  }+\frac{\left(
n^{2}-2m\right)  }{2}\left(  n-2\right)  +4mn-n^{3}\\
&  =-\frac{n^{2}}{2}\frac{n-2}{n-1}-\frac{2\left(  n-2\right)  m}{\left(
n-1\right)  }+\frac{2m^{2}}{n-1}+\left(  n-2\right)  m\\
&  =\frac{n\left(  n-2\right)  }{\left(  n-1\right)  }\left(  2m-\frac{n^{2}%
}{2}\right)  +\frac{2m^{2}}{n-1}+\left(  n-2\right)  m<\frac{2m^{2}}%
{n-1}+\left(  n-2\right)  m.
\end{align*}

This completes the proof of (\ref{bo4}) and of Theorem \ref{th2}.
\end{proof}

\textbf{Acknowledgement} The author is grateful to Cecil Rousseau for his
patient explanation of Olpp's paper and for interesting discussions.

\end{document}